\documentclass{amsart}
\usepackage{amsmath}
\usepackage{amssymb}
\begin{document}

\title[Triangular Equilibrium Points in the Generalized ... --Drag]{Triangular Equilibrium Points in  the Generalized  Photogravitational Restricted Three Body Problem   with Poynting-Robertson drag }\footnotetext[1]{\textit{Review Bull.Cal. Math. Soc., 12 (1 \& 2) 109-114 (2004)}}
\author[B.S. Kushvah and B. Ishwar]{ B.S.Kushvah$^1$ and B.Ishwar$^2$ \\ 
 (1)JRF Under DST Project (2) P. I. DST Project\\
 B.R.A. Bihar University Muzaffarpur-842001, India \\
email: bskush@hotmail.com}
\maketitle

\begin{abstract}
  In this paper,we have found the equations of motion of Generalized Photogravitational restricted three body problem with Poynting-Robertson drag. The problem is generalized in the sense that smaller primary is supposed to be an oblate spheroid. The bigger primary is considered as radiating. The equations of motion are affected by radiation pressure force oblateness and P-R drag.We have located triangular equilibrium points in our problem. All classical results involving photogravitational and oblateness in restricted three body problem may be verified from this result.\\
\end{abstract}
\noindent
{\bf AMS Classification:70F15}\\
{\bf Keywords}:Triangular Points/Generalized Photogravitational/RTBP/P-R drag.

\newcommand{\abc}{{(1-\frac{5}{2}A_2) }}
\newcommand{\aac}{{(1-\frac{A_2}{2}) }}
\newcommand{\amc}{{(1-\mu) }}
\newcommand{\adc}{{\frac{\delta ^2}{2} }}
\section{Introduction} 
%{\bf {1. Introduction}}\\
 In general three body problem, we study the motion of three finite bodies.The problem is restricted in the sense that one of the three masses is taken to be so small that the gravitational effect on the other masses by third mass is negligible. The small body is known as infinitesimal mass and remaining two massive bodies as finite masses or primaries. The classical restricted three body problem is generalized to include the force of radiation pressure,the Poynting-Robertson drag effect and oblateness effect.\\

Poynting(1903) has stated that the particle such as small meteors or cosmic dust  are comparably affected by gravitational and light  radiation force, as they approach luminous celestial bodies. He also suggested that infinitesimal body in solar orbit suffers a gradual loss of angular momentum and ultimately spiral into the Sun. In a system of coordinates where the Sun is at rest, radiation scattered by infinitesimal mass in the direction of motion suffers a blue shift and in the opposite direction it is red shifted. This gives rise to net drag force which opposes the direction of motion. The proper relativistic treatment of this problem  was formulated by Robertson(1937) who showed that to first order in ${\vec{\frac{V}{c}}}$ the radiation pressure force is given by\\
\begin{equation}\label{eq:F}
\vec{F}=\displaystyle{{F_p}\biggl\{ \frac{\vec{R}}{R} - \frac{\vec{V}.\vec{R}\vec{R}}{cR^{2}}- \frac{\vec{V}}{c} \biggr\}}
\end{equation}

where $F_{p}$=$\frac{3Lm}{16{\pi}R^{2}\rho{sc}}$   denotes the measure of the radiation pressure force,$\vec{R}$ the position vector of $P$ with respect to radiation sours $S$, $\vec{V}$ the corresponding velocity vector and $c$ the velocity of light. In the  expression of $F_{p}$, $L$ is luminosity of the radiating body, while $m, \rho $ and $s$ are the mass, density and cross section of the particle respectively.

The first term in equation (~\ref{eq:F}) express the radiation pressure. The second term represents the Doppler shift of the incident radiation and the third term is due to the absorption and subsequent re-emission of part of  the incident radiation. These last two terms taken together are the Poynting-Robertson effect. The Poyinting-Robertson effect will operate to sweep small particles of the solar system into the Sun at cosmically rapid rate.

Stanley P., Wyatt J. R. and Fred L. Whipple (1950) have shown that the P-R effect has been of very little significance galacticaly.Its importance is restricted to small particles orbiting in the vicinity of individual stars. The importance of the radiation influence on celestial bodies has been recognized by many scientists such as Kozai (1961), McCraken and Alexander (1968), Ferrz-Mello (1972), Simmons {\it et al.}(1985), Vkrouhlixky (1993, 1994), Murray C.D.(1994), Ragos O.and Zafiropoulos F.A. (1995).\\

 Colombo {\it et al}. (1996) analyzed the stability of the equilibrium points in the presence of radiation pressure which include the Poynting-Robertson drag term. They showed that the points were unstable to such a drag force.Chernikove (1970) has dealt with the Sun-Planet-Particle model. He concludes that the P-R effect renders unstable those liberation points known to be conditionally stable in the classical case.  Schuerman D.W. (1980) has studied the triangular points. He has shown that these points are unstable on time scale long compared the period of revolution of the two massive bodies.\\

 This problem has an interesting application for artificial satellite and future space colonization. It has been suggested that the classical triangular points of the Sun-Jupiter or Sun-Earth system would be convenient sites to locate future space colonies.

In this paper we consider Sun-Planet-Particle model with Sun as a radiating body, planet as an oblate spheroid. We have found that the coordinates of equilibrium points are the functions of mass reduction factor  $ q_{1} $, P-R drag $W_{1}$ and coefficients of oblateness $A_{2}$. All classical results in restricted three body problem involving radiation and oblateness may be deduced from this result.\\

\section{Equations of Motion}
%{\bf 2. Equations of Motion:}\\

\- We suppose $m_{1}$, the mass of more massive radiating primary and $m_{2}$ the mass of smaller primary which is an oblate spheroid. Let these primaries revolve around their centre of mass in circular orbits in the plane of motion which coincides to the equatorial plane of $m_{2}$.\\

\-  We consider the barycentric rotating co-ordinate system$OXYZ$ relative to inertial system with angular velocity $\omega$ and common $Z$ -axis.We have taken line joining the primaries as $X$ -axis. $OX$ and $OY$ in the equatorial plane of $m_{2}$ and $OZ$ coinciding with the polar axis of $m_{2}$. Let $r_{e}$ and$ r_{p}$ be the equatorial and polar radii of $m_{2}$  and $r$ be the distance between primaries.Let infinitesimal mass $m$ be placed at a point $P(X,Z,Y)$.Then potential at $P$ due to $m_{1}$ and  $m_{2}$ is 
\begin{equation}\label{eq:V}
V=-\frac{k^{2}mm_{1}q_{1}}{r_{1}}-\frac{k^{2}mm_{2}}{r_{2}}-\frac{k^{2}mm_{2}A_{2}}{2r^{2}_{3}}\end{equation}
where $k^{2}$is the Gaussian constant of gravitational and $A_{2}=\frac{r^{2}_{e}-r^{2}_{p}}{5r^{2}}$,oblateness coefficient of $m_{2}$, $r_{1}$ and $r_{2}$ are the distances from $m_{1}$ to $P$ and $m_{2}$ to$P$ respectively. We take units such that sum of the masses and distance between primaries  as unity. The unit of time i.e. time period of $m_{1}$ about $m_{2}$ will consists of $2\pi$ units such that $k=1$. Then perturbed mean motion of the primaries is $n^{2}=1+\frac{3A_{2}}{2}$. Let $\mu=\frac{m_{2}}{m_{1}+m_{2}}$ then $\frac{m_{1}}{m_{1}+m_{2}}=1-\mu$ with $m_{1}>m_{2}$ where $\mu$ is mass parameter.Then (~\ref{eq:F}) and (~\ref{eq:V}) and use the same technique as Chernikove (1970) and Schuerman (1980). In dimensionless coordinate system, the dimensionless velocity of light will be given by $c_{d}=c$, which depends on the physical masses of the primaries and the distance between them.

For inertial reference system the total acceleration acting on $P$ is given by Coriolis relation 
\begin{eqnarray}
&&\vec{a}+2\vec{\omega}\times\vec{v}+\vec{\omega}\times(\vec{\omega\times{r}})=-\displaystyle{\frac{(1-\mu)\vec{r_1}}{r^3_1}-\frac{\mu{\vec{r_2}}}{r^3_2}-\frac{3}{2}\frac{\mu{A_2}\vec{r_2}}{r^5_2}}\notag\\
&&+\displaystyle{\frac{(1-\mu)(1-q_1)}{r^2_1}\biggl\{\frac{\vec{r_1}}{r_1}-\frac{(\vec{\dot{r_1}}+\vec{\omega}\times\vec{r_1}).\vec{r_1}\vec{r_1}}{c_d r^2_1}-\frac{\vec{\dot{r_1}}+\vec{\omega}\times{\vec{r_1}}}{c_d}\biggr\}}\notag
\end{eqnarray}
where
\begin{gather*} 
\vec{r_{1}}=x\hat{i}+y\hat{j},\quad
\vec{v}=\dot{x}\hat{i}+\dot{y}\hat{j},\quad
\vec{a}=\ddot{x}\hat{i}+\ddot{y}\hat{j},\quad
\vec{\omega}=n\hat{k},\quad
\vec{r_{1}}=(x+\mu)\hat{i}+y\hat{j},\\
\vec{r_{2}}=(x+\mu-1)\hat{i}+y\hat{j},\quad
r^{2}_{1}=(x+\mu)^{2}+y^{2},\quad
r^{2}_{2}=(x+\mu-1)^{2}+y^{2}.  
\end{gather*}
 Substituting all these values in above relation and comparing the components of $\hat{i}$ and $\hat{j}$,we get the equations of motion of the infinitesimal mass particle in $xy$-plane.\\

\begin{eqnarray}
U_x=\ddot{x}-2n\dot{y}&=&n^{2}x-\frac{(1-\mu)q_1(x+\mu)}{r^3_1}-\frac{\mu(x+\mu-1)}{r^3_2}-\frac{3}{2}\frac{\mu{A_2}(x+\mu-1)}{r^5_2} \notag\\
&&-\frac{W_1}{r^2_1}\biggl\{\frac{(x+\mu)}{r^2_1}[(x+\mu){\dot{x}+y\dot{y}}] +\dot{x}-ny \biggr\}\label{eq:ux}\\\notag
U_y=\ddot{y}+2n\dot{x}&=& n^{2}y-\frac{(1-\mu)q_1 y}{r^3_1}-\frac{\mu{y}}{r^3_2}-\frac{3}{2}\frac{\mu{A_2}y}{r^5_2} \notag\\&&-\frac{W_1}{r^2_1}\biggl\{\frac{y}{r^2_1}[(x+\mu)\dot{x}+y\dot{y}]+\dot{y}+n(x+\mu)\biggr\}\label{eq:uy}
\end{eqnarray}
where $W_1=\frac{(1-\mu)(1-q_1)}{c_d}$,$n^2=1+\frac{3}{2}A_2$,$q=1-\frac{F_p}{F_g}$ is a mass reduction factor expressed in terms of the particle radius $a$,density $\rho$ radiation pressure efficiency factor $\chi$ (in C.G.S. system) $q=1-\frac{5.6\times{10^{-5}}}{a\rho}\chi$. The assumption $q$=constant is equivalent to neglecting fluctuations in the beam of solar radiation and the effect of the planets shadow. Obviously $q\leq1$. Equations of motion (~\ref{eq:ux}) and (~\ref{eq:uy}) can be written as
\begin{align*}
\ddot{x}-2n\dot{y}=\frac{\partial{U_1}}{\partial{x}}+F_x ,\quad \ddot{y}+2n\dot{x}=\frac{\partial{U_1}}{\partial{y}}+F_y 
\end{align*}
\- where 
\begin{equation*}
U_1=\frac{n^2}{2}{(x^2+y^2)}+\displaystyle{\frac{(1-\mu){q_1}}{r_1}}+\frac{\mu}{r_2}+\displaystyle{\frac{\mu{A_2}}{2{r_2^3}}}
\end{equation*}
\begin{align}
F_x&=-\displaystyle{\frac{W_1}{r^2_1}\biggl\{\frac{(x+\mu)}{r^2_1}[(x+\mu)\dot{x}+y\dot{y}]+\dot{x}-ny\biggr\}}\label{eq:fx} \\[1ex]
\text{and}\quad 
F_y&=-\displaystyle{\frac{W_1}{r^2_1}\biggl\{\frac{y}{r^2_1}[(x+\mu)\dot{x}+y\dot{y}]+\dot{y}+n(x+\mu)\biggr\}}\label{eq:fy}
\end{align}
$F_x$,$F_y$ are the partial derivatives of drag function with respective to $x$ and $y$ respectively, which are purely functions of the particle's position and velocity.

Now multiplying equations (~\ref{eq:fx}) by $2\dot{x}$,(~\ref{eq:fy}) by $2\dot{y}$ and adding, we get,
\[
2\dot{x}\ddot{x}+2\dot{y}\ddot{y}=2\biggl(\dot{x}\frac{\partial{U_1}}{\partial{x}}+\dot{y}\frac{\partial{U_1}}{\partial{y}}\biggr)+2\bigl(\dot{x}F_x+\dot{y}F_y\bigr)
\] 
which can be written as $\frac{dC}{dt}=-2\bigl(\dot{x}F_{x}+\dot{y}F_{y}\bigr)$, where $C=2U_1-{\dot{x}}^2-{\dot{y}}^2$. The quantity $C$ is Jacobi Integral. The zero velocity curves are given by $C=2U_{1}(x,y)$

\section{ Location of Triangular Equilibrium Points}
 
 For the triangular equilibrium points $y\neq{0}$, $U_x=U_y=0$ then from equations (~\ref{eq:ux}) and (~\ref{eq:uy})
\begin{eqnarray}
n^{2}x-\displaystyle{\frac{(1-\mu){q_1}(x+\mu)}{r^3_1}-\frac{\mu(x+\mu-1)}{r^3_2}-\frac{3}{2}\frac{\mu{A_2}(x+\mu-1)}{r^5_2}+\frac{W_1}{r^2_1}ny}&=&0\label{eq:ux0}\\
n^{2}y-\displaystyle{\frac{(1-\mu){q_1}y}{r^3_1}-\frac{\mu{y}}{r^3_2}-\frac{3}{2}\frac{\mu{A_2}y}{r^5_2}-\frac{W_1}{r_1}n(x+\mu)}&=&0\label{eq:uy0}
\end{eqnarray}
Multiplying equations (~\ref{eq:ux0}) by $y$, (~\ref{eq:uy0}) by $(x+\mu)$ then subtracting, we get
\begin{equation*}
\biggl\{n^2-\displaystyle{\frac{1}{r^3_2}\bigl(1+\frac{3A_2}{2r^2_2}\bigr)}\biggr\}\mu{y}=nW_1.
\end{equation*}
In the case of photogravitational restricted three body problem we have $r_{1_0}=q^{1/3}=\delta$ (say) and $r_{2_0}=q_2^{1/3}=1$.So we suppose that due to P-R drag and oblateness, perturbation in $r_{1_0}$ and $r_{2_0}$ are $\epsilon_1$ and $\epsilon_2$ respectively where ${\epsilon_{i}}'$s are very small,

\begin{equation}r_1=q^{1/3}_1(1+\epsilon_1) \quad \text{and} \quad r_2=1+\epsilon_2\end{equation}

Putting these values in equation (~\ref{eq:ux0}) and (~\ref{eq:uy0}) and neglecting higher order terms of small quantities, we get 
\[
\epsilon_2{y}=\displaystyle{\frac{nW_1}{3\mu(1+\frac{5}{2}A_2)}}
\]
 we may write as $\epsilon_2{y}=\epsilon_2{y_0}$,i.e. $\epsilon_2=\displaystyle{\frac{nW_1(1-\frac{5}{2}A_2)}{3\mu{y_0}}}$,$\epsilon_1=-\displaystyle{\frac{nW_1}{6(1-\mu){y_0}}}-\frac{A_2}{2}$
 Hence, we get
 \begin{align*}
r_1&=\delta\biggl\{1-\displaystyle{\frac{nW_1}{6(1-\mu){y_0}}}-\frac{A_2}{2}\biggr\},r_2=1+\displaystyle{\frac{nW_1}{3\mu{y_0}}(1-\frac{5}{2}A_2)}
 \end{align*}
Since  $x+\mu=\frac{r^2_1-r^2_2+1}{2}$ and $y^2=r^2_1-(x+\mu)^2$,

\begin{eqnarray}
&&x=x_0\biggl\{1-\displaystyle{\frac{nW_1[\amc\abc+\mu\aac\adc]}{3\mu\amc y_0 x_0}}-\adc\frac{A_2}{x_0}\biggr\}\label{eq:L4x}\\
&&y=y_0\biggl\{1-\displaystyle{\frac{nW_1\delta^2[2\mu-1-\mu(1-\frac{3A_2}{2})\adc+7\amc\frac{A_2}{2}]}{3\mu\amc y^3_0}}-\displaystyle{\frac{\delta^2(1-\adc)A_2}{y^2_0}}\biggr\}^{1/2}\label{eq:L4y}
\end{eqnarray}
$(x_0,y_0)$ are coordinates of $L_4,L_5$ in the photogravitational restricted three body problem where \[
x_0=\adc-\mu,\quad y_0=\pm\delta\biggl(1-\frac{\delta^4}{4}\biggr)^{1/2}, \quad \delta=q^{1/3}_1
\]
Equations (~\ref{eq:L4x}) and (~\ref{eq:L4y}) are valid for $W_1\ll 1,A_2\ll 1$.

\section{ Conclusion}
%{\bf 4. Conclusion :}
 The equilibrium points $L_4$, and $L_5$ are given by equations (~\ref{eq:L4x}) and (~\ref{eq:L4y}).
{\bf Case 1.}
   We see that when P-R effect not included i.e.,$W_1=0$ then we get $x=x_0-\adc A_2$ and $y=y_0\{1-\adc(1-\adc)\frac{A_2}{y^2_0}\}$ which are the $L_4$ and $L_5$ points in the case of generalized photogravitational restricted three body problem.
{\bf Case 2.}
When $A_2=0$ then $n=1$, i.e., small primary is spherically symmetric. We have 
\[
x=x_0\biggl\{1-\displaystyle{\frac{W_1[\amc+\mu\adc]}{3\mu\amc x_0 y_0}}\biggr\},\quad y=y_0\biggl\{1-\displaystyle{\frac{W_1\delta^2[2\mu-1-\mu\adc]}{6\mu\amc y^3_0}}\biggr\}.
\]
this result coincides with {\it Schuerman(1980)}.
{\bf Case 3.} When  $A_2=0,W_1=0$,then we have for $ q_1=1=\delta,$ $x=\frac{1}{2}-\mu$, $y=\pm\frac{\sqrt{3}}{2}$
These are the coordinates of classical restricted three body problem.

Finally we conclude that position of triangular equilibrium points are affected by mass reduction factor, P-R drag and oblateness coefficient.

{\bf Acknowledgment:}\\

We are thankful to D.S.T. Government of India New Delhi for sanctioning a project DST/MS/140/2K on this topic.

\end{document}